\documentclass[a4paper, 11pt]{article}
\usepackage{stmaryrd}

\usepackage{amscd}
\usepackage{amssymb,amsmath,amsthm}
\usepackage{mathptmx}
\usepackage{mathrsfs}
\usepackage{color}

\DeclareSymbolFont{letters}{OML}{cmm}{m}{it}

\DeclareMathAlphabet{\mathcal}{OMS}{cmsy}{m}{n}

\newtheorem{theorem}{Theorem}[section]
\newtheorem{corollary}[theorem]{Corollary}
\newtheorem{lemma}[theorem]{Lemma}

\newtheorem{remark}[theorem]{\textrm{\textbf{Remark}}}

\theoremstyle{definition}
\newtheorem{definition}[theorem]{Definition}

\definecolor{darkgreen}{rgb}{0,.7,0}

\begin{document}

\title{Expressing Second-order Sentences in Intuitionistic Dependence Logic\footnote{This work is part of the European Science Foundation EUROCORES LogICCC project Logic for Interaction (LINT). The research leading to the current article was supported by the Finnish National Graduate School in Mathematics and its Applications.}}
\author{Fan Yang\\
Department of Mathematics and Statistics \\
University of Helsinki, Finland\\
fan.yang@helsinki.fi}

\date{}

\maketitle

\begin{abstract}
Intuitionistic dependence logic was introduced by Abramsky and V\"{a}\"{a}n\"{a}nen \cite{AbVan09} as a variant of dependence logic under a general construction of Hodges' (trump) team semantics.  It was proven that there is a translation from intuitionistic dependence logic sentences into second order logic sentences. In this paper, we prove that the other direction is also true, therefore intuitionistic dependence logic is equivalent to second order logic on the level of sentences.\end{abstract}


\section{Introduction}

Dependence Logic (\textbf{D}), as a new approach to independence
friendly logic (IF-logic) \cite{HintikkaSandu1989}, was introduced
in \cite{Van07dl}. Hodges gave a compositional semantics for
IF-logic in \cite{Hodges1997a}, \cite{Hodges1997b}, trump semantics (or team semantics). Recent
research by Abramsky, V\"{a}\"{a}n\"{a}nen
\cite{AbVan09} generalized Hodges' construction for 
team semantics (or trump semantics) and introduced BID-logic, which extends dependence
logic and includes both intuitionistic implication and linear
implication, as well as intuitionistic disjunction. We call the intuitionistic fragment of BID-logic
``intuitionistic dependence logic (\textbf{ID})''. In this paper,
we study the expressive power of sentences of intuitionistic dependence logic. By the method of \cite{Enderton1970} and \cite{Walkoe1970}, we know that sentences of \textbf{D} have exactly the same expressive power as sentences of $\Sigma^1_1$, the existential second order fragment.  
It was proven in \cite{AbVan09} that \textbf{ID} sentences are
expressible in second-order logic (\textbf{SO}). We will show that the other direction is also true; that is, there is a translation from the sentences of the full \textbf{SO} into \textbf{ID}. In particular, \textbf{D} sentences (or $\Sigma^1_1$ sentences) are expressible in \textbf{ID}. This means that \textbf{ID} is so powerful that it is equivalent to the full \textbf{SO} on the level of sentences.

We name the logic under discussion ``intuitionistic'' dependence logic, because the implication of this logic satisfies the axioms of the usual intuitionistic implication. As we know, the usual intuitionistic logic (either propositional or first-order) is weaker than classical logic (see e.g. \cite{vanDalen01}). However, although the idea of introducing the intuitionistic implication in the general context of Hodges' construction is very natural, as the above-mentioned result of this paper shows, in the team semantics context, \textbf{ID} is stronger than (classical) \textbf{D}. It is worthwhile to point out that restricted to (classical) first-order
formulas (or flat formulas), \textbf{ID} is in fact classical. It is only between rednon-classical formulas (dependence
formulas) that the intuitionistic implication does play a role.

Another logic in the team semantics setting with the same expressive power as the full second order logic is the so-called \emph{team logic}, which is the logic of dependence logic extended with classical negation (see \cite{Van07dl} and also \cite{Ville_thesis}, \cite{KontinenNurmi2009}). The significance of intuitionistic dependence logic is that the equivalence of \textbf{ID} and the full \textbf{SO} on the sentence level is established without the presence of the logical connective classical negation.

Throughout the paper, we assume readers are familiar with the
standard Tarskian semantics of first-order logic and the standard
semantics of second-order logic. We assume that the domain of a first-order model $M$ is non-empty. For any model $M$, an \emph{assignment} $s$ on $M$ is a function from a finite set $dom(s)$ of variables into the domain of $M$ ($dom(s)$ will be always clear from the context). Let $a$ be an element in $M$, and $x$ a variable in $dom(s)$. We write $s(a/x)$ for the assignment with $dom(s(a/x))=dom(s)\cup\{x\}$ which agrees with $s$ everywhere except that it maps $x$ to $a$. The sequences of variables $\langle x_{i,1},\dots,x_{i,n}\rangle$ and $\langle x^i_1,\dots,x^i_n\rangle$ with subscripts and superscripts are abbreviated as $\overline{x_i}$ and $\overline{x^i}$, respectively. Similarly for sequences of constants and elements of
models. For any assignment $s$ for $\overline{x}$, we write $s(\overline{x})$ for the sequence
$\langle s(x_1),\dots,s(x_n)\rangle$. We use the standard
abbreviation $\forall\overline{x}$ to stand for a sequence of universal
quantifiers $\forall x_1\dots\forall x_n$ (the length of $\overline{x}$
is always clear from the context or does not matter); similarly for
existential quantifiers.

\section{Intuitionistic Dependence Logic as a Fragment of BID-logic}

In this section, we define intuitionistic dependence logic, which
is the intuitionistic fragment of BID-logic introduced in
\cite{AbVan09}. We will also recall some basic properties of dependence
logic proved in \cite{Van07dl} in
this general framework.

BID-logic is obtained from a general construction of Hodges' team
semantics (or trump semantics) \cite{Hodges1997a}, \cite{Hodges1997b}. Well-formed
formulas of BID-logic (in negation normal form) are given by the following grammar
\[\begin{split}
    \phi::= & \,\alpha\mid=\!\!(t_1,\dots,t_n)\mid
\neg=\!\!(t_1,\dots,t_n)\mid\bot\mid\phi\wedge\phi\mid
\phi\otimes\phi\mid\phi\ovee\phi \mid\\
      &  \phi\to\phi\mid \phi\multimap\phi\mid\forall x\phi\mid\exists
x\phi
  \end{split}
\] 
where $\alpha$ is a first-order atomic or negated atomic formula (first-order literal), $t_1,\dots,t_n$
are terms. Formulas of the form $=\!\!(t_1,\dots,t_n)$ are called dependence atomic formulas or dependence atoms. The
disjunctions ``$\otimes$'' and ``$\ovee$'' are called \emph{split disjunction} and \emph{intuitionistic disjunction}, the implications ``$\to$'' and ``$\multimap$'' are called
\emph{intuitionistic implication} and \emph{linear implication}, respectively.
 The set $Fv(\phi)$ of free variables of a formula $\phi$ of BID-logic is defined in the standard way except for the case of dependence atoms:
\[Fv(=\!\!(t_1,\dots,t_n))=Var(t_1)\cup\dots\cup Var(t_n),\]
where $Var(t_i)$ ($1\leq i\leq n$) is the set of variables occurring in $t_i$.
We call $\phi$ a sentence in case $Fv(\phi) = \emptyset$.

For the semantics for BID-logic, we adopt and generalize
Hodges' team semantics. For any model $M$, a \emph{team} $X$ of
$M$ is a set of assignments on $M$ with the same domain $dom(X)$. We define two
operations on teams. For any team $X$ of $M$, and any function
$F:X\to M$, the \emph{supplement team} $X(F/x)=\{s(F(s)/x):s\in
X\}$ and the \emph{duplicate team} $X(M/x)=\{s(a/x): a\in M,
~s\in X\}$. Now we give the team semantics for BID-logic. For any suitable model $M$ and any team $X$ of $M$ whose domain includes the set of free variables of the formula under discussion,
\begin{itemize}
\item $M\models_X\alpha$ with $\alpha$ first-order literal iff
for all $s\in X$, $M\models_s\alpha$ in the usual Tarskian
  semantics sense;
  \item $M\models_X=\!\!(t_1,\dots,t_n)$ iff for all $s,s'\in X$ with $s(t_1)=s'(t_1)$,$\dots$, $s(t_{n-1})=s'(t_{n-1})$, it holds that
  $s(t_n)=s'(t_n)$;
  \item $M\models_X\neg =\!\!(t_1,\dots,t_n)$ iff $X=\emptyset$;
  \item $M\models_X\bot$ iff $X=\emptyset$;
  \item $M\models_X\phi\wedge\psi$ iff $M\models_X\phi$ and
  $M\models_X\psi$;
  \item $M\models_X\phi\otimes\psi$ iff there exist teams $Y,Z\subseteq X$ with $X=Y\cup Z$ such that $M\models_Y\phi$ and
  $M\models_Z\psi$;
  \item $M\models_X \phi\ovee\psi$ iff $M\models_X \phi$ or $M\models_X
    \psi$;
  \item $M\models_X \phi\to\psi$ iff for any team $Y\subseteq X$, if $M\models_Y \phi$ then $M\models_Y\psi$;
\item $M\models_X \phi\multimap\psi$ iff for any team $Y$ with $dom(Y)=dom(X)$, if $M\models_Y \phi$ then $M\models_{X\cup Y}\psi$;
  \item $M\models_X\exists x\phi$ iff $M\models_{X(F/x)}\phi$ for some function $F:X\to
  M$;
  \item $M\models_X\forall x\phi$ iff $M\models_{X(M/x)}\phi$.
\end{itemize}
We say that a formula $\phi$ is \emph{satisfied} by a team $X$ in a model $M$, if $M\models_X\phi$ holds. It can be easily shown that satisfaction of a formula of BID-logic depends only on the interpretations of the variables occurring free in the formula. Sentences have no free variable and there is only one assignment with empty domain, namely the empty assignment $\emptyset$. We say that a sentence $\phi$ is \emph{true} in $M$ if the team $\{\emptyset\}$ of empty assignment satisfies $\phi$, i.e. $M\models_{\{\emptyset\}}\phi$.  We use the standard notation $\phi\models \psi$ to mean that for any suitable model $M$, $M\models \phi$ implies $M\models \psi$.

The intuitionistic implication and linear implication are adjoints
of the corresponding conjunctions; that is
\[\phi\wedge\psi\models\chi\Longleftrightarrow \phi\models \psi\to\chi,\]
\[\phi\otimes\psi\models\chi\Longleftrightarrow \phi\models
\psi\multimap\chi.\]
 The propositional variant of BID-logic without
dependence formulas is essentially the BI logic, 
the ``logic of
Bunched Implications'' introduced in \cite{OHearnPym1999},
\cite{Pym2002}. The fragment with connectives $\wedge$, $\otimes$
and quantifiers is the usual \emph{dependence logic}, where as pointed out in \cite{AbVan09}, $\otimes$ is in fact the multiplicative conjunction instead of disjunction (therefore the notation ``$\otimes$''). The
intuitionistic fragment of BID-logic is called \emph{intuitionistic
dependence logic}. More precisely,  well-formed formulas of
\textbf{D} are formed by the following grammar
\[    \phi::= \alpha\mid=\!\!(t_1,\dots,t_n)\mid\neg=\!\!(t_1,\dots,t_n)\mid\phi\wedge\phi\mid\phi\otimes\phi \mid \forall
x\phi\mid\exists x\phi
\]
where $\alpha$ is a first-order literal and $t_1,\dots,t_n$ are
terms. Well-formed formulas of \textbf{ID} are formed by the
following grammar
\[    \phi::= \alpha\mid=\!\!(t)\mid
\bot\mid\phi\wedge\phi\mid\phi\ovee\phi \mid \phi\to\phi\mid\forall
x\phi\mid\exists x\phi
\]
where $\alpha$ is a first-order atomic formula and $t$ is a term. Note that
the dependence atoms of \textbf{ID} have only single terms and $M\models_X=\!\!(t)$ means intuitively that $t$ behaves in team $X$ as a constant. We will see later in Lemma \ref{eq-BID} that dependence atom with several variables $=\!\!(t_1,\dots,t_n)$ is definable by the constancy dependence atoms $=\!\!(t_i)$.

The most important property of BID-logic is the \emph{downwards
closure} property that for any formula $\phi$, if $M\models_X\phi$
and $Y\subseteq X$, then $M\models_Y\phi$. A formula $\phi$ is said
to be \emph{flat} if for all suitable models $M$ and teams $X$
\[M\models_X\phi\Longleftrightarrow (M\models_{\{s\}}\phi\mbox{ for all }s\in X).\]

We call
the \textbf{D} formulas with no occurrence of dependence
subformulas \emph{(classical) first-order formulas} (of BID-logic), that is first-order formulas of BID-logic are formulas with only first-order literals, $\wedge$, $\otimes$, $\forall x$ and $\exists x$. Throughout the paper, we sometimes talk about first-order formulas of BID-logic and the usual first-order formulas at the same time, in such cases, we identify the first-order connective $\otimes$ of BID-logic with the usual first-order connective $\vee$.

\begin{lemma}\label{flatness} First-order formulas are flat.
\end{lemma}
\begin{proof}
Easy, by induction on the structure of formulas.
\end{proof}

\begin{lemma}\label{sent-flat} Sentences of BID-logic are flat.
\end{lemma}
\begin{proof} To evaluate sentences with no free variables, we only consider the singleton team $\{\emptyset\}$ of the empty assignment $\emptyset$. 
\end{proof}

It is easy to observe that \textbf{D} and \textbf{ID} have the \emph{empty team property}, that is for any \textbf{D} or \textbf{ID} formula $\phi$, any model $M$, the empty team satisfies $\phi$, i.e. $M\models_\emptyset \phi$. However, the full BID-logic does not have such property. For example, for any model $M$, $M\not\models_\emptyset (x=x)\multimap (x\neq x)$.

\section{First-order Formulas are Expressible in
Intuitionistic Dependence Logic}

In this section, we show that every first-order formula (of BID-logic) is
logically equivalent to a formula in intuitionistic dependence
logic. Two formulas $\phi$
and $\psi$ of BID-logic are said to be \emph{logically equivalent} to each other, in symbols $\phi\equiv\psi$, if for any suitable model $M$ and team $X$ with $dom(X)\supseteq Fv(\phi)\cup Fv(\psi)$, it holds that
\[M\models_X\phi\Longleftrightarrow M\models_X\psi.\]

\begin{lemma}\label{eq-BID}
We have the following logical equivalences in BID-logic:
\begin{description}
  \item[(1)] $=\!\!(t_1,\dots,t_n) \,\equiv\, =\!\!(t_1)\wedge \dots\wedge =\!\!(t_{n-1})\to
  =\!\!(t_n)$ for any terms $t_1,\dots,t_n$;
  \item[(2)] $\neg\phi \,\equiv\,\phi\to
  \bot$ whenever $\phi$ is an atom (first-order or dependence atom);
  \item[(3)] $(\phi\to\bot)\to\bot\,\equiv\,\phi$ whenever $\phi$ is a flat formula;
  \item[(4)] $\phi\otimes\psi\,\equiv\,(\phi\to\bot)\to\psi$ whenever both $\phi$ and $\psi$ are
flat formulas.
\end{description}
\end{lemma}
\begin{proof}
Items (1)-(3) can be proved easily. We only show item (4). That is to
show that for any model $M$ and any team $X$ with $dom(X)\supseteq
Fv(\phi)\cup Fv(\psi)$ it holds that
\[M\models_X\phi\otimes\psi\Longleftrightarrow M\models_X(\phi\to\bot)\to\psi.\]

$\Longrightarrow$: Suppose $M\models_X\phi\otimes\psi$. Then there
exist two teams $Y,Z$ with $X=Y\cup Z$ such that $M\models_Y\phi$ and
$M\models_Z\psi$. For any nonempty $U\subseteq X$ with
$M\models_U\phi\to\bot$, downwards closure gives that for any $s\in
U$, $M\models_{\{s\}}\phi\to\bot$, i.e. $M\not\models_{\{s\}}\phi$.
Since $M\models_Y\phi$, in view of the downwards closure we
conclude that $s\not\in Y$, thus $U\subseteq Z$, which implies
$M\models_U\psi$ by downwards closure.

$\Longleftarrow$: Suppose $M\models_X(\phi\to\bot)\to\psi$. Define
\begin{center}
$Y=\{s\in X\mid M\models_{\{s\}} \phi\}$ and $Z=\{s\in X\mid
M\not\models_{\{s\}} \phi\}$.
\end{center}
Clearly, $X=Y\cup Z$. If $Z\neq \emptyset$, then for any $s\in Z\subseteq X$, we have that
$M\models_{\{s\}}\phi\to\bot$, thus since $M\models_{\{s\}}
(\phi\to\bot)\to\psi$, we obtain that $M\models_{\{s\}}\psi$. Now
both $M\models_Y\phi$ and $M\models_Z\psi$ follow from the flatness
of $\phi$ and $\psi$.
\end{proof}

\begin{remark}\label{IDL-classical-singleton}
Item (2) of the above lemma shows that for atomic formulas, dependence negation and intuitionistic negation have the same meaning. Moreover, restricted to singleton teams, the connectives $\to$ and $\ovee$ of \textbf{ID}
behave as classical connectives. To determine whether a flat formula is satisfied by a team $X$, one in fact only needs to consider the satisfaction of singleton subteams of $X$, for which the intuitionistic negation behaves classically. This explains why items (3), (4) hold.
\end{remark}

Next we define expressibility.
\begin{definition}
Let $\mathscr{L}$ be a sublogic of BID-logic. We say that a formula
$\phi$ of BID-logic is \emph{expressible} in $\mathscr{L}$, if there
exists an $\mathscr{L}$ formula $\psi$ such that $\phi\equiv\psi$.
\end{definition}

\begin{theorem}\label{FO-IDL}First-order formulas (of BID-logic) are expressible in
\textbf{ID}.
\end{theorem}
\begin{proof}Assuming that every first-order formula is in prenex normal form and the negation-free part is in conjunctive normal form, the theorem follows immediately from items (2), (4) of Lemma
\ref{eq-BID}. 
 
For example, the
quantifier-free first-order formula $(\alpha\otimes(\neg\beta\otimes\gamma))\wedge
\delta$ in conjunctive normal form, where
$\alpha,\beta,\gamma,\delta$ are first-order atoms, can be
translated into \textbf{ID} as follows:
\[\begin{split}
    ( \alpha\otimes(\neg\beta\otimes\gamma))\wedge \delta & \Longrightarrow ((\alpha\to \bot)\to (\neg\beta\otimes
\gamma))\wedge\delta\\
&\Longrightarrow  ((\alpha\to \bot)\to ((\neg\beta\to\bot)\to
\gamma))\wedge\delta \\
      & \Longrightarrow ((\alpha\to \bot)\to (((\beta\to\bot)\to\bot)\to
\gamma))\wedge\delta.
 \end{split}
\]
\end{proof}

\section{Sentences of Dependence Logic and $\Sigma^1_1$ are
Expressible in Intuitionistic Dependence Logic}

It follows from \cite{Enderton1970} and \cite{Walkoe1970} that IF-logic is equivalent to the $\Sigma^1_1$ fragment of
\textbf{SO} on the level of sentences. Dependence logic, which is equivalent to
IF-logic on the level of sentences is therefore also equivalent to the $\Sigma^1_1$ fragment.
V\"{a}\"{a}n\"{a}nen in \cite{Van07dl} gave a direct translation from
one logic into the other. In this section, we will prove that there exists a translation from sentences of \textbf{D}  or $\Sigma^1_1$ sentences into \textbf{ID}. This proof will be generalized in Section 5 to the full second-order logic.

\begin{definition}
Let $\mathscr{L}_{\textbf{SO}}$ and $\mathscr{L}_{\rm{BID}}$ be sublogics of
second-order logic and of BID-logic, respectively.
\begin{enumerate}
\item We say that a sentence
$\phi$ of $\mathscr{L}_{\rm{BID}}$ is \emph{expressible} in
$\mathscr{L}_{\textbf{SO}}$, if there exists an $\mathscr{L}_{\textbf{SO}}$ sentence
$\psi$ such that for any suitable model $M$,
\[M\models\psi\Longleftrightarrow M\models_{\{\emptyset\}}\phi.\]
\item  We say that a sentence
$\psi$ of $\mathscr{L}_{\textbf{SO}}$ is \emph{expressible} in
$\mathscr{L}_{\rm{BID}}$, if there exists an $\mathscr{L}_{\rm{BID}}$ sentence
$\phi$ such that for any suitable model $M$,
\[M\models\psi\Longleftrightarrow M\models_{\{\emptyset\}}\phi.\]
\end{enumerate}
\end{definition}

\begin{theorem}[\cite{Van07dl}]
\textbf{D} sentences are expressible in the $\Sigma^1_1$ fragment of
\textbf{SO}.
\end{theorem}

 We sketch the proof of the next theorem. In the next
section, we will generalize the translation in the next theorem to translate all second
order sentences, first into BID-logic, and in the end into
\textbf{ID}.

\begin{theorem}[\cite{Van07dl}]\label{Sigma11-DL} $\Sigma^1_1$ sentences are expressible in
\textbf{D}.
\end{theorem}
\begin{proof}(idea)
Without loss of generality, we may assume every $\Sigma_1^1$
sentence $\phi$ is of the following special Skolem normal form
\[\exists f_1\dots\exists f_n\forall x_1\dots\forall x_m\psi,\]
where $\psi$ is quantifier-free and first-order, and for each $1\leq i\leq n$, every occurrence of the function symbol $f_i$ is of the same form
$f_i x_{i_1}\dots x_{i_q}$ for some fixed sequence $\langle x_{i_1},\dots, x_{i_q}\rangle$ of variables from the set $\{x_1,\dots,x_m\}$. We find a \textbf{D} sentence $\phi^\ast$ which expresses $\phi$. The idea behind the sentence $\phi^\ast$ is that in $\phi$, we replace each occurrence of the function symbol $f_i$ by a new variable $y_i$, and add a dependence atom to specify that $y_i$ is functionally determined by the arguments $x_{i_1},\dots, x_{i_q}$ of $f_i$. This can be done because we have required that each occurrence of $f_i$ is of the same form $f_i x_{i_1}\dots x_{i_q}$. To be precise, the \textbf{D} sentence $\phi^\ast$ is defined as follows:
\begin{equation}\label{DL-NF}
\begin{split}
    \phi^\ast\,:=\,\forall x_1\dots \forall x_m\exists y_1\dots \exists y_n & (\,=\!\!(x_{1_1},\dots, x_{1_q},y_1)\wedge \\
      & \dots \wedge=\!\!(x_{n_1},\dots,x_{n_q},y_n)\wedge\psi'\,),
  \end{split}
\end{equation}
 where $\psi'$
is obtained from $\psi$ by replacing everywhere $f_ix_{i_1}\dots
x_{i_q}$ by the new variable $y_i$ for each $1\leq i\leq n$. In $\phi^\ast$, the dependence atoms together with the existential quantifiers enable us to pick exactly those functions corresponding to the functions assigned to the existentially quantified function variables $f_1,\dots,f_n$ in $\phi$.\end{proof}

\begin{remark}
Equation (\ref{DL-NF}) with the first-order quantifier-free formula $\psi'$ in conjunctive normal form is a
normal form for \textbf{D} sentences.
\end{remark}

Note that in the normal form (\ref{DL-NF}) of a \textbf{D} sentence, the only
subformulas that are not in the language of \textbf{ID} are
dependence atoms with several variables and first-order quantifier-free formulas. As we have proved in the
previous section, these two kinds of formulas are both expressible
in \textbf{ID}. Therefore we obtain the next theorem.

\begin{theorem}\label{DL-s-IDL}
\textbf{D} sentences are expressible in \textbf{ID}.
\end{theorem}
\begin{proof}
Let $\phi$ be a \textbf{D} sentence in the normal form (\ref{DL-NF}). The \textbf{ID}
sentence that expresses $\phi$ is obtained by replacing the subformulas of the
form $=\!\!(x_{i,1},\dots,x_{i,q},y_i)$ by the formula
$=\!\!(x_{i_1})\wedge \dots\wedge =\!\!(x_{i,q})\to
  =\!\!(y_i)$ and the first-order quantifier-free formula $\psi'$ by its equivalent \textbf{ID} formula obtained from Theorem
  \ref{FO-IDL}.
\end{proof}

\begin{theorem}\label{Sigma11-IDL}
$\Sigma_1^1$ sentences are expressible in \textbf{ID}.
\end{theorem}
\begin{proof}
Follows from Theorem \ref{Sigma11-DL} and Theorem \ref{DL-s-IDL}.
\end{proof}

Negation of dependence logic, as well as that of BID-logic do not satisfy
Law of Excluded Middle (neither for split disjunction nor for intuitionistic disjunction) and are therefore not the classical negation. 
Dependence logic extended with classical negation (denoted by ``$\thicksim$'') is called \emph{team
logic}. Definable team properties of team logic correspond exactly to all second order properties, in particular, sentences of team logic have the same expressive power as sentences of the full second order logic, see \cite{Van07dl}, \cite{Ville_thesis}, \cite{KontinenNurmi2009} for further discussions
on team logic. The result of the equivalence of \textbf{ID} and the full \textbf{SO} sentences to be presented in the next section shows that \textbf{ID}, as a logic without classical negation, is actually equivalent to team logic on the level of sentences. However, team logic and \textbf{ID} are not equivalent in general, since for example, the classical negation of $\bot$, denoted by $\thicksim\!\bot$, is a well-formed formula in team logic, while in \textbf{ID}, it is not expressible. Indeed, suppose the classical negation of $\bot$ was expressible by an \textbf{ID} formula $\phi$, namely for any suitable model $M$ and team $X$,
\[M\not\models_X \bot\Longleftrightarrow M\models_X \phi.\]
In particular, for empty team $\emptyset$, it would hold that 
\[M\not\models_\emptyset \bot\Longleftrightarrow M\models_\emptyset \phi,\]
but this is never the case since by the semantics $M\models_\emptyset \bot$ and by empty set property of \textbf{ID}, $M\models_\emptyset \phi$ always holds. This also shows that classical negation in general is not definable in \textbf{ID}. However, for sentences of \textbf{ID}, as shown in the next lemma, the intuitionistic negation does give a certain kind of ``classical'' negation.

\begin{lemma}\label{IDL-negation}
For any sentence $\phi$ of BID-logic, we have that for any suitable model $M$
\[M\models_{\{\emptyset\}}\phi\to \bot\Longleftrightarrow M\not\models_{\{\emptyset\}}\phi.\]
\end{lemma}
\begin{proof}
Trivial.
\end{proof}

A sentence of BID-logic has no free variables, thus it is said to be true if and only if the team $\{\emptyset\}$ of the empty assignment (a singleton team) satisfies it. As pointed out in Remark
\ref{IDL-classical-singleton}, when
restricted to singleton teams, the semantics of \textbf{ID} is in
fact classical. This explains why such kind of ``classical'' negation is definable in \textbf{ID}.

Using the intuitionistic (``classical'') negation for sentences, we
are able to express $\Pi^1_1$ sentences as well.

\begin{corollary}\label{Pi11-IDL}$\Pi^1_1$ sentences
are expressible in \textbf{ID}.
\end{corollary}
\begin{proof}
Let $\psi$ be a $\Pi^1_1$ sentence. Note that $\psi$ is equivalent
to $\neg\phi$ for some  $\Sigma^1_1$ sentence $\phi$. By
Theorem \ref{Sigma11-IDL}, there exists an \textbf{ID} sentence
$\phi^\ast$ such that
\[M\models\phi\Longleftrightarrow
M\models_{\{\emptyset\}}\phi^\ast\] for all suitable models $M$. Since
$\phi^\ast$ is a sentence, by Lemma \ref{IDL-negation}, we have
that
\[M\models_{\{\emptyset\}}\phi^\ast\to\bot\Longleftrightarrow M\not\models_{\{\emptyset\}}\phi^\ast\Longleftrightarrow M\not\models\phi\Longleftrightarrow M\models\neg\phi\Longleftrightarrow M\models\psi,\]
thus $\phi^\ast\to \bot$ is the sentence of \textbf{ID}
expressing $\psi$.
\end{proof}

\section{Second-order Sentences are Expressible in
Intuitionistic Dependence Logic}

In this section, we will generalize the proofs of Theorem
\ref{Sigma11-DL} and Theorem \ref{Sigma11-IDL} to show that sentences of the full \textbf{SO} are expressible in \textbf{ID}.
Together with the result of the next theorem, proved in
\cite{AbVan09}, we will be able to conclude that the
expressive power of \textbf{ID} sentences is so strong that it is, in fact, 
equivalent to that of sentences of the full \textbf{SO}.

\begin{theorem}[\cite{AbVan09}]\label{IDL-SO}
\textbf{ID} sentences are expressible in second-order logic.
\end{theorem}

In order to proceed to the main theorem of this paper (Theorem
\ref{SOL-IDL}), we first recall the normal form of \textbf{SO}
formulas. 

\begin{theorem}\label{SO-NF}[Normal Form of \textbf{SO}]\
Every second order sentence is equivalent to a formula of the form
\[\forall\overline{f^1}\exists\overline{f^2}\dots\forall \overline{f^{2n-1}}\exists\overline{f^{2n}}\forall\overline{x}\psi,\]
where $\psi$ is quantifier-free, and we assume without loss of generality that for the corresponding $Q\in\{\forall,\exists\}$, each $Q\overline{f^i}=Qf^i_1\dots Qf^i_p$ and each $f^i_j$ is of arity $q$.
\end{theorem}

The basic idea of the translation for sentences of the full \textbf{SO} is generalized from that of the proof of Theorem \ref{Sigma11-DL} for $\Sigma^1_1$ sentences. For each \textbf{SO} sentence in a special normal form (to be clarified in Lemma \ref{niceform_lm}), we replace each function variable by a new variable and specify the functionality of the new variable by adding the corresponding dependence atoms. We have seen in the proof of Theorem \ref{Sigma11-DL} that dependence atoms together with existential quantifiers enable us to simulate existentially quantified function variables; on the other hand, universally quantified function variables can also be simulated using dependence atoms and intuitionistic implications. In this way, we will be able to express all \textbf{SO} sentences in \textbf{ID}.

To make this idea work, we need to first turn every \textbf{SO} sentence $\phi$ into a better normal form than the one in Theorem \ref{SO-NF}, that is we need to guarantee that for each $q$-ary function variable $f^i_j$, every occurrence of $f^i_j$ in $\phi$ is of the same form $f^i_jx_{i,j_1}\dots x_{i,j_q}$ for some fixed sequence $\langle x_{i,j_1}\dots x_{i,j_q}\rangle$ of variables (this normal form is inspired by the $\Sigma^1_1$ normal form in Theorem \ref{Sigma11-DL}, see Section 6.3 in \cite{Van07dl} for detailed discussions). To this end, we need three lemmas.

The first lemma removes nesting of function symbols in a formula.
\begin{lemma}\label{niceform_func}Let $\phi(ft_1\dots t_q)$ be any first-order formula, in which the $q$-ary function symbol $f$ has an occurrence of the form $ft_1\dots t_q$ for some terms $t_1\dots t_q$. Then we have that
\[\models\phi(ft_1\dots t_q)\leftrightarrow\forall x_1\dots\forall x_q((t_1=x_1)\wedge\dots\wedge(t_q=x_q)\to\phi(fx_1\dots x_q)),\]
where $x_1,\dots,x_q$ are new variables and $\phi(fx_1\dots x_q)$
is the formula obtained from $\phi(ft_1\dots t_q)$ by replacing
everywhere $ft_1\dots t_q$ by $fx_1\dots x_q$.
\end{lemma}
\begin{proof}Easy.
\end{proof}


The second lemma unifies the arguments of function symbols in a formula.
\begin{lemma}\label{niceform2_func}Let $\phi(fx_1\dots x_q,fy_1\dots y_q)$ be a first-order formula, in which the $q$-ary function symbol $f$ has an occurrence of the form $fx_1\dots x_q$ and an occurrence of the form $fy_1\dots y_q$ with $\{x_1\dots x_q\}\cap\{y_1\dots y_q\}=\emptyset$. Then we have that
\[\begin{split}
    \models & \forall x_1\dots \forall x_q\forall y_1\dots\forall y_q\phi(fx_1\dots
x_q,fy_1\dots y_q) \\
      & \leftrightarrow \exists g\forall x_1\dots \forall
x_q\forall y_1\dots\forall y_q(\phi(fx_1\dots x_q,gy_1\dots y_q)\\
&\wedge((x_1=y_1)\wedge\dots\wedge(x_q=y_q)\to
(fx_1\dots x_q=gy_1\dots y_q))),
  \end{split}
\] where
$\phi(fx_1\dots x_q,gy_1\dots y_q)$ is the first-order formula obtained from the formula $\phi(fx_1\dots,x_q,fy_1\dots y_q)$ by replacing everywhere $fy_1\dots y_q$ by $gy_1\dots
y_q$.
\end{lemma}
\begin{proof}Easy.
\end{proof}

The next lemma gives a nice normal form for \textbf{SO} sentences. 
\begin{lemma}\label{niceform_lm}
Every \textbf{SO} formula is equivalent to a formula $\phi$ of the
form
\[\forall f^1_1\dots \forall f^1_{p}\exists f^2_1\dots\exists f^2_{p}\dots\,\dots \forall f^{2n-1}_1\dots \forall f^{2n-1}_{p}\exists f^{2n}_1\dots \exists f^{2n}_{p}\forall x_1\dots\forall x_m\psi,\]
where
\begin{itemize}
  \item $\psi$ is quantifier free;
  \item for each $1\leq i\leq 2n$ and each $1\leq j\leq p$, every occurrence of the $q$-ary function symbol $f^i_j$  is of the same form $f^i_j\mathbf{x}^{i,j}$, where
  \[\mathbf{x}^{i,j}=\langle x_{i,j_1},\dots ,x_{i,j_{q}}\rangle\]
 with $\{x_{i,j_1},\dots ,x_{i,j_{q}}\}\subseteq \{x_1,\dots,x_m\}$.
\end{itemize}
\end{lemma}
\begin{proof}
By applying Lemma \ref{niceform_func}  and Lemma \ref{niceform2_func} several times and adding dummy quantifiers to the \textbf{SO} formulas in the normal form described in Theorem \ref{SO-NF}. 
\end{proof}



The next lemma states that under the right valuations, the behavior of functions can be simulated by new variables. This technical lemma will play a role in the proof of Lemma
\ref{Pi2n-IDL-lm}.

\begin{lemma}\label{func-team}
Let $\phi(\overline{f},\overline{x})$ be any first-order quantifier-free formula with function symbols $f_1,\dots,f_p$,
where every occurrences of each $q$-ary function symbol $f_j$ is of the same form
\[f_jx_{j_1}\dots x_{j_q},\]
where $\{ x_{j_1}\dots x_{j_q}\}\subseteq \{x_1,\dots,x_m\}$. Let $(M,\overline{F})$ be any suitable model with new function symbols $\underline{f_1},\dots,\underline{f_p}$ interpreted as $F_1,\dots,F_p$, respectively. Let $y_1,\dots,y_p$ be new variables and $s$ an assignment with domain 
\[\{x_1,\dots,x_m,y_1,\dots,y_p\}\] 
such that for all $1\leq j\leq p$,
\begin{equation}\label{func-team-eq1}
s(y_j)=F_j(s(x_{j_1}),\dots,s(x_{j_q})).
\end{equation}
Then
\[(M,\overline{F},s(\overline{x}))\models \phi(\overline{f},\overline{x})\Longleftrightarrow M\models_{\{s\}}\phi',\]
where $\phi'$ is the first-order formula of BID-logic obtained from $\phi$ by
replacing everywhere $f_jx_{j_1}\dots x_{j_q}$ by $y_j$ for
each $1\leq j\leq p$.
\end{lemma}
\begin{proof}
It is easy to show by induction that for any term $t$,
$s(t)=s(t')$, where $t'$ is obtained from $t$ by replacing everywhere $f_jx_{j_1}\dots x_{j_q}$ by $y_j$ for
each $1\leq j\leq p$. Next, we show the lemma by induction on $\phi$.
The only interesting case is the case that $\phi\equiv
\psi\vee\chi$. In this case, we have that
\begin{align*}
    (M,\overline{F},s(\overline{x}))\models \psi\vee\chi \Longleftrightarrow&
    (M,\overline{F},s(\overline{x}))\models \psi\mbox{ or }(M,\overline{F},s(\overline{x}))\models \chi\\
     \Longleftrightarrow& M\models_{\{s\}}\psi'\mbox{ or }M\models_{\{s\}} \chi'\\
    &\mbox{(by induction hypothesis)}\\
    \Longleftrightarrow& M\models_{\{s\}}\psi'\otimes\chi'\footnotemark\\
    &(\text{since }\{s\}=\{s\}\cup \{s\}=\{s\}\cup\emptyset).
\end{align*}
\footnotetext{Here, as mentioned in Section 2, we identify the connective disjunction ``$\vee$'' in \textbf{SO} and the split disjunction ``$\otimes$'' in BID-logic.}
\end{proof}

Now we are in a position to give the translation from \textbf{SO} sentences into \textbf{ID}. To simplify notations, this will be done in two steps. In the first step (Lemma \ref{Pi2n-IDL-lm}), for each \textbf{SO} formula $\phi$, we find an equivalent sentence $\phi^\ast$ in BID-logic. The second step (Theorem \ref{Pi2n-IDL}) will turn $\phi^\ast$ finally into an equivalent \textbf{ID} sentence.

If $\mathbf{x}=\langle x_1,\dots,x_n\rangle$ is a sequence of variables, then we abbreviate the atomic formula $=\!\!(x_1,\dots,x_n,y)$ as
$=\!\!(\mathbf{x},y)$. If $X$ is a team of $M$, then the duplicate team $X(M/x_1)\dots(M/x_n)$ is abbreviated as $X(M/x_1,\dots,x_n)$.

\begin{lemma}\label{Pi2n-IDL-lm}
\textbf{SO} sentences are expressible in BID-logic.
\end{lemma}
\begin{proof}
Without loss of generality, we may assume that every \textbf{SO}
sentence $\phi$ is of the form described in Lemma \ref{niceform_lm}.
For each pair $\langle i,j\rangle$ ($1\leq i\leq 2n$, $1\leq j\leq p$), pick a new variable $u_{i,j}$ not occurring in $\phi$. We inductively define BID formulas $\delta_i$ for $2n\geq i\geq 1$ as follows: 
\[\delta_{2n}~:=~ \exists u_{2n,1}\dots\exists u_{2n,p}(\Theta_{2n}\wedge \psi'),\]
for $2n>i\geq 1$,  
\[\delta_{i}~:=~ \left\{
                         \begin{array}{ll}
                            \Theta_{i}\to \delta_{i+1}, & \hbox{if $i$ is odd;} \\
                           \exists u_{i,1}\dots\exists u_{i,p}(\Theta_{i}\wedge \delta_{i+1}), & \hbox{if $i$ is even,}
                         \end{array}
                       \right.
\]
where
\[\Theta_i=\bigwedge_{j=1}^{p}\,=\!\!(\mathbf{x}^{i,j},u_{i,j})\]
 and $\psi'$ is the first-order
formula obtained from $\psi$ by replacing everywhere each $f^i_j\mathbf{x}^{i,j}$ by $u_{i,j}$. Let
\[\phi^\ast=\forall u_{1,1}\dots \forall u_{1,p}\forall u_{3,1}\dots \forall u_{3,p}\dots\,\dots \forall u_{2n-1,1}\dots \forall u_{2n-1,p}\forall \overline{x}\delta_{1}\]
\begin{align*}
\big[\text{i.e. }\phi^\ast=&\forall u_{1,1}\dots \forall u_{1,p}\forall u_{3,1}\dots \forall u_{3,p}\dots\,\dots \forall u_{2n-1,1}\dots \forall u_{2n-1,p}\forall \overline{x}\\
&(\Theta_1\to \exists u_{2,1}\dots\exists u_{2,p}(\Theta_2\wedge(\Theta_3\to \exists u_{4,1}\dots\exists u_{4,p}(\Theta_4\wedge\cdots\,\cdots\\
&~~~~\cdots\,\cdots \wedge(\Theta_{2n-1}\to \exists
u_{2n,1}\dots\exists
u_{2n,p}(\Theta_{2n}\wedge\psi'\underbrace{))\cdots\,\cdots))))}_{2n}\big].
 \end{align*}

The general idea behind the BID formula $\phi^\ast$ is that the $\delta_i$'s for $i$ odd, simulate the $\forall \overline{f^i}$'s, and the $\delta_i$'s for $i$ even, simulate the $\exists\overline{f^i}$'s in the \textbf{SO} sentence $\phi$. The rest of the proof is devoted to show that such sentence
$\phi^\ast$ does express the \textbf{SO} sentence $\phi$,
i.e. to show that for any suitable model $M$ it holds that
\begin{equation*}
M\models\phi\Longleftrightarrow M\models_{\{\emptyset\}}\phi^\ast.
\end{equation*}

\ \\

``$\Longrightarrow$'': Suppose $M\models\phi$. Then for any sequence of
functions 
\[F^1_1,\dots,F^1_{p}:M^{q}\to M,\]
there
exists a sequence of functions (depending on $\overline{F^1}$)
\[F^2_1{\scriptstyle(\,\overline{F^1}\,)},\dots,F^2_{p}{\scriptstyle(\,\overline{F^1}\,)}: M^{q}\to M\]
such that for any 
$\dots$ $\dots$
for any sequence of functions
\[F^{2n-1}_1,\dots,F^{2n-1}_{p}:M^q\to M,\]
there
exists a sequence of functions (depending on $\overline{F^1},\overline{F^3}\dots,\overline{F^{2n-1}}$\,)
\[F^{2n}_1{\scriptstyle(\,\overline{F^1},\dots,\overline{F^{2n-1}})},\dots,F^{2n}_{p}{\scriptstyle(\overline{F^1},\dots,\overline{F^{2n-1}}\,)}: M^q\to M\]
such that
\begin{equation}\label{Pi12n-1}
(M,\overline{F^1},\dots,\overline{F^{2n}})\models
\forall
\overline{x}\psi(\,\overline{f^1},\dots,\overline{f^{2n}}\,).
\end{equation}

Let $Y_1$ be a nonempty subteam of
\[X=\{\emptyset\}(M/\overline{u_{1}},\overline{u_{3}},\dots, \overline{u_{2n-1}},\overline{x})\] such
that $M\models_{Y_1}\Theta_1$. 
It suffices to show that
\begin{equation}\label{Pi12n-0}
 M\models_{Y_1} \delta_2\text{, i.e. }   M\models_{Y_1} \exists u_{2,1}\dots\exists u_{2,q}(\Theta_2\wedge\delta_{3}).
 \end{equation}


The team $Y_1$ corresponds to a sequence of functions $F^1_1{\scriptstyle(Y_1)},\dots,F^1_p{\scriptstyle(Y_1)}:M^q\to M$ defined as follows: for any $1\leq j\leq p$, and for some fixed $a_0\in M$, let
\[F^1_j(\overline{d})=\left\{
                         \begin{array}{ll}
                           s(u_{1,j}), & \hbox{if there exists $s\in Y_1$ such that $s(\mathbf{x}^{1,j})=\overline{d}$;} \\
                           a_0\in M, & \hbox{otherwise.}
                         \end{array}
                       \right.
\]
Each $F^1_j$ is well-defined. Indeed, for any $\overline{d}\in M^q$, any $s,s'\in Y_1$ such that
\[s(\mathbf{x}^{1,j})=\overline{d}=s'(\mathbf{x}^{1,j}),\]
since $M\models_{Y_1}\, =\!\!(\mathbf{x}^{1,j},u_{1,j})$, we must have that
\[s(u_{1,j})=s'(u_{1,j}).\]

Now, using the functions $F^2_1{\scriptstyle( \overline{F^1})},\dots,F^2_p{\scriptstyle( \overline{F^1})}$, we define a sequence of functions $\alpha_{2,1}{\scriptstyle(F^2_1)},\dots,\alpha_{2,p}{\scriptstyle(F^2_p)}$ from the corresponding supplement teams of $Y_1$ to $M$ such that the supplement team $Y_1 (\alpha_{2,1}/u_{2,1})\dots(\alpha_{2,p}/u_{2,p})$ satisfies $\Theta_2\wedge \delta_3$. For each $1\leq j\leq p$, define the function
\[    \alpha_{2,j}:Y_1 (\alpha_{2,1}/u_{2,1})\dots
(\alpha_{2,j-1}/u_{2,j-1})\to M\]
corresponding to $F^2_j{\scriptstyle(\overline{F^1})}$ by taking
\[\alpha_{2,j}(s)=F^2_j(s(\mathbf{x}^{2,j})).\]
Put
\[Y_2=Y_1(\alpha_{2,1}/u_{2,1})\dots
(\alpha_{2,p}/u_{2,p}).\] It suffices to show that
$M\models_{Y_2}\Theta_2$ and
\begin{equation}\label{Pi12n-00}
M\models_{Y_2}\delta_3\text{, i.e. } M\models_{Y_2}\Theta_3\to \delta_{4}.
\end{equation}

The former is obvious by the definitions of $Y_2$ and
$\overline{\alpha_{2}}$. To show the latter, repeat the same argument and construction $n-1$ times, and it then suffices to show that for any nonempty subteams $Y_3$ of $Y_2$, $Y_5$ of $Y_4$, $\dots$ , $Y_{2n-1}$ of
$Y_{2n-2}$ such that 
\begin{center}
$M\models_{Y_3}\Theta_3$, $M\models_{Y_5}\Theta_5$, $\dots$,
$M\models_{Y_{2n-1}}\Theta_{2n-1}$, 
\end{center}
it holds that
\begin{equation}\label{Pi12n-2}
M\models_{Y_4}\Theta_4, M\models_{Y_6}\Theta_6,\dots,M\models_{Y_{2n}}\Theta_{2n}
\end{equation}
 and
$M\models_{Y_{2n}}\psi'$, where $Y_4,Y_6\dots,Y_{2n}$ are supplement teams defined in the
same way as above. Clause (\ref{Pi12n-2}) follows immediately from the definitions of $Y_4,Y_6,\dots,Y_{2n}$ and
$\overline{\alpha_{4}},\overline{\alpha_{6}},\dots
\overline{\alpha_{2n}}$.
To show $M\models_{Y_{2n}}\psi'$, since $\psi'$ is flat (first-order), it suffices
to show $M\models_{\{s\}}\psi'$ hold for all $s\in Y_{2n}$. 

For the functions $\overline{F^1}{\scriptstyle(Y_1)},\overline{F^2}{\scriptstyle(\overline{F^1})},\dots,\overline{F^{2n-1}}{\scriptstyle (Y_{2n-1})},\overline{F^{2n}}{\scriptstyle (\overline{F^{2n-1}})}$ obtained as above,
by (\ref{Pi12n-1}) we have that
\[
    (M,\overline{F^1},\dots,\overline{F^{2n}}, s(\overline{x}))\models
\psi(\overline{f^1},\dots,\overline{f^{2n}},\overline{x}).
\]
Now, it follows from the definitions of $\overline{F^1},\dots,\overline{F^{2n}}$ that condition (\ref{func-team-eq1}) in Lemma \ref{func-team}  is satisfied for each $F^i_j$, hence,  an application of Lemma \ref{func-team} gives the desired result that $M\models_{\{s\}}\psi'$.

\ \\

``$\Longleftarrow$'': Suppose $M\models_{\{\emptyset\}}\phi^\ast$. Then
\[   M\models_X  \Theta_1\to \delta_{2},\]
where 
\[X=\{\emptyset\}(M/\overline{u_{1}},\overline{u_{3}},\dots, \overline{u_{2n-1}},\overline{x}).\]
Let $F^1_1,\dots,F^1_{p}:M^{q}\to M$ be an arbitrary sequence of functions. Take a subteam $Y_1{\scriptstyle (\overline{F^1})}$ of $X$ which corresponds to these functions by putting
\[ \begin{split}Y_1=  \{s\in \{\emptyset\}(M/\overline{u_{1}},&\overline{u_3}\dots\overline{u_{2n-1}},\overline{x})\\
& \mid s(u_{1,1})=F^1_1(s(\mathbf{x}^{1,1})),\dots,s(u_{1,p})=F^1_p(s(\mathbf{x}^{1,p}))\}.   \end{split} \] 

Clearly, $M\models_{Y_1}\Theta_1$ holds,
thus we have that $M\models_{Y_1}\delta_2$ holds (i.e., (\ref{Pi12n-0}) holds).
So there exist functions
\[\alpha_{2,1}{\scriptstyle (\overline{F^1})}:Y_1\to M,\dots\,\dots,\alpha_{2,p}{\scriptstyle (\overline{F^1})}:Y_1(\alpha_{2,1}/u_{2,1})\dots
(\alpha_{2,p-1}/u_{2,p-1})\to M\]
depending on $\overline{F^1}$ such that $M\models_{Y_2}\Theta_2$ and $M\models_{Y_2}\delta_3$ holds (i.e., (\ref{Pi12n-00}) holds),
 where 
 \[Y_2=Y_1(\alpha_{2,1}/u_{2,1})\dots(\alpha_{2,p}/u_{2,p}).\]
Now, we define functions $F^2_1{\scriptstyle (\overline{F^1})},\dots, F^2_{p}{\scriptstyle (\overline{F^1})}:M^q\to M$, which simulate $\alpha_{2,1},\dots, \alpha_{2,p}$ as follows: for each $1\leq j\leq p$ and for any $\overline{d}\in M^q$, let
\[F^2_j(\overline{d})=s(u_{2,j})\text{ for some }s\in Y_2\text{ such that }s(\mathbf{x}^{2,j})=\overline{d}.\]
Note that the definition of $Y_2$ guarantees such $s$ in the above definition always exists, and moreover, each $F^2_j$ is well-defined since for any $s,s'\in Y_2$ with
\[s(\mathbf{x}^{2,j})=\overline{d}=s'(\mathbf{x}^{2,j}),\]
as $M\models_{Y_2}=\!\!(\mathbf{x}^{2,j},u_{2,j})$, we must have that
\[s(u_{2,j})=s'(u_{2,j}).\]

Repeat the same argument and construction $n-1$ times to define inductively for any sequences of
functions $\overline{F^3},\overline{F^5},\dots,\overline{F^{2n-1}}$, the subteams $Y_3$
of $Y_2$, $\dots$ , $Y_{2n}$ of $Y_{2n-1}$ such that
\[M\models_{Y_3}\Theta_3,~M\models_{Y_5}\Theta_5,\dots,M\models_{Y_{2n-1}}\Theta_{2n-1},\]
and the supplement teams $Y_4,Y_6,\dots,Y_{2n}$ satisfy
\[M\models_{Y_4}\Theta_4,~M\models_{Y_6}\Theta_6,\dots,M\models_{Y_{2n-2}}\Theta_{2n-2},~M\models_{Y_{2n}}\Theta_{2n}\wedge \psi',\]
and to define inductively 
the sequences of functions
\[\overline{F^4},\overline{F^6}\dots,\overline{F^{2n}}:M^{q}\to M,\]
according to the functions $\overline{\alpha_4},\overline{\alpha_6},\dots,\overline{\alpha_{2n}}$ obtained from the existential quantifiers $\exists \overline{u_{4}},\exists \overline{u_6}\dots,\exists \overline{u_{2n}}$. It then suffices to show that 
\[(M,\overline{F^1},\dots,\overline{F^{2n}})\models\forall \overline{x} \psi(\overline{f^1},\dots,\overline{f^{2n}}).\]

Let $\overline{a}$ be an arbitrary sequence in $M$ of the same length as that of $\overline{x}$. By the construction of $Y_{2n}$, there must exists $s\in Y_{2n}$ such that $s(\overline{x})=\overline{a}$. Since $M\models_{Y_{2n}}\psi'$, by downwards closure, we have that $M\models_{\{s\}}\psi'$. Note that by the definitions of $\overline{F^1},\dots,\overline{F^{2n}}$, condition (\ref{func-team-eq1}) in Lemma \ref{func-team}  is satisfied for each $F^i_j$, hence, an application of Lemma \ref{func-team} gives the desired result that 
\[(M,\overline{F^1},\dots,\overline{F^{2n}},s(\overline{x}))\models\psi(\overline{f^1},\dots,\overline{f^{2n}},\overline{x}).\]
\end{proof}

Observe that in the sentence $\phi^\ast$ in the proof of Theorem
\ref{Pi2n-IDL-lm}, the only subformulas that are not in the
language of \textbf{ID} are dependence atoms with several variables and first-order
formulas, both of which are expressible in \textbf{ID}. Therefore
we are able to turn it into an equivalent \textbf{ID} sentence.

\begin{theorem}\label{Pi2n-IDL}
\textbf{SO} sentences are expressible in \textbf{ID}.
\end{theorem}
\begin{proof}
Follows from Lemma \ref{Pi2n-IDL-lm}, Lemma \ref{eq-BID} and
Theorem \ref{FO-IDL}.
\end{proof}

Finally, we arrive at the following theorem.

\begin{theorem}\label{SOL-IDL}
\textbf{SO} sentences are expressible in \textbf{ID}, and vice versa.
\end{theorem}
\begin{proof}
Follows from Theorem \ref{Pi2n-IDL} and \ref{IDL-SO}.
\end{proof}

\begin{remark}
It is easy to observe that in the above translation, the intuitionistic disjunction $\ovee$ did not play a role. In fact, $\ovee$ is definable in \textbf{ID} uniformly by the expression $\phi\ovee\psi\equiv \theta_1\wedge\theta_2$,
where
\[\theta_1:=\forall x\forall y(x= y)\to((\phi\to\bot)\to\psi)\]
(which deals with the case that the model has only one element) and
\[
\begin{split}\theta_2:=\forall x\exists y(x=y\to\bot)\to \exists w\exists u\big(=\!\!(w)\wedge=\!\!(u)&\wedge \left((w=u)
\to\phi\right)\\
&\wedge\left((w= u\to\bot)\to\psi\right)\big)
\end{split}
\]
(which deals with the other cases).
\end{remark}

\begin{remark}
In fact, Lemma \ref{niceform_lm} gives a normal form for $\Pi^1_{2n}$ ($n\in\omega$) sentences, therefore the \textbf{ID} sentences of the form in Theorem \ref{Pi2n-IDL} can be viewed as \textbf{ID} normal form for $\Pi^1_{2n}$ sentences. Moreover, using the ``classical negation'' of \textbf{ID} for sentences, applying the same trick as that in the proof of Corollary \ref{Pi11-IDL}, one can obtain an \textbf{ID} normal form for $\Sigma^1_{2n}$ sentences $\phi$, namely $\psi^\ast\to\bot$, for $\psi^\ast$ the translation of the $\Pi^1_{2n}$ sentence $\psi\equiv \neg\phi$. 

Using the same trick as that in Lemma \ref{niceform_lm}, one can also obtain a nice normal form for $\Sigma^1_{2n-1}$ sentences. Thus, the above observation (normal form) holds for $\Sigma^1_{2n-1}$ and $\Pi^1_{2n-1}$ sentences as well. In particular, the proof of Theorem \ref{Sigma11-IDL} (for $\Sigma^1_1$ sentences) can then be viewed as a special case of the proof of Theorem \ref{Pi2n-IDL}.
\end{remark}

\section{Further Work}

\subsection{Expressive Power of \textbf{ID} Open Formulas}
In this paper, we proved that the expressive power of sentences of \textbf{ID} is equivalent to that of sentences of the full second order logic. In \cite{KontVan09}, it was proven that with respect to nonempty teams, open formulas of \textbf{D} defines exactly those properties that are definable in $\Sigma^1_1$ with an extra predicate, occurring only negatively, for the nonempty teams. A similar result can be obtained for open formulas of \textbf{ID}, that is with respect to nonempty teams, open formulas of \textbf{ID} defines exactly those properties that are definable in the full second order logic with an extra predicate, occurring only negatively,  for the nonempty teams, see \cite{Yang2010} for details. 

\subsection{Linear Dependence Logic}
One other interesting fragment of BID-logic is the linear dependence logic (\textbf{LD}). This is the fragment extended from dependence logic by adding the linear implication. More precisely, well-formed formulas of
\textbf{LD} are formed by the following rule
\[    \phi::= \alpha\mid=\!\!(t_1,\dots,t_n)\mid\neg=\!\!(t_1,\dots,t_n)\mid\phi\wedge\phi\mid\phi\otimes\phi \mid \phi\multimap\phi\mid \forall
x\phi\mid\exists x\phi
\]
where $\alpha$ is a first-order literal and $t_1,\dots,t_n$ are
terms. One may wonder whether the translation discussed in this paper works for \textbf{LD} sentences. This is still unclear, however, one has to realize in the first place that for \textbf{LD} sentences, the situation is more complicated. Because, as pointed out in Section 2, \textbf{LD} does not have the empty team property. For any sentence $\phi$ of BID-logic, we define the truth value $\llbracket\phi\rrbracket$ of $\phi$ on a model $M$ to be 
\[\llbracket\phi\rrbracket=\{X\mid M\models_X\phi,~X\in\wp (M^{Fv(\phi)})\}.\]
The failure of empty team property implies that for any \textbf{LD} sentence $\phi$, the truth value $\llbracket\phi\rrbracket$ lies in the three-element set $\{\,\{\{\emptyset\},\emptyset\},\{\emptyset\},\emptyset\,\}$, namely, \textbf{LD} is three-valued (see \cite{AbVan09} for more details).

One result we have obtained along this line is that the similar translation does apply to the second order $\Pi^1_2$ fragment and \textbf{LD} in the following sense. Let $\phi$ be a $\Pi^1_2$ sentence of the form described in Lemma \ref{niceform_lm} and $\phi^\ast$ the sentence of BID-logic defined in Lemma \ref{Pi2n-IDL-lm} which expresses $\phi$. We replace the intuitionistic implication $\to$ in $\phi^\ast$ by the linear implication $\multimap$ and denote the resulting \textbf{LD} sentence by $\phi^{\ast\ast}$. Then for any model $M$, the following is true:
\[M\models \phi\,\Longleftrightarrow\, M\models_\emptyset \phi^{\ast\ast}.\]
This shows that sentences of \textbf{LD} goes beyond $\Sigma^1_1$ already with respect to the truth value $\emptyset$ . 

\subsection{IF-logic and the Full Second Order Logic}
It is well-known that Independence friendly logic is equivalent to $\Sigma^1_1$, thus to \textbf{D}, on the level of sentences. This indicates a possibility of obtaining the same result of this paper for an extension of IF-logic. However, the result of this paper relies heavily on the role the intuitionistic implication plays in the translation; that is, it is based on a deep understanding of the general framework of Hodges' team semantics. Since the original semantics of IF-logic was given by means of imperfect information games (\cite{HintikkaSandu1989}), to obtain the same result of this paper for a reasonable extension of IF-logic, one may have to seek for a different notion, a game-theoretic one, which corresponds to the intuitionistic implication in the team semantics context.

\paragraph{Acknowledgements.} 

The author would like to thank Juha Kontinen and Jouko V\"{a}\"{a}n\"{a}nen for very helpful discussions on the topic of this paper.


\begin{thebibliography}{14}


\bibitem{AbVan09}
\textsc{Abramsky, S.}, and \textsc{J. V\"{a}\"{a}n\"{a}nen}, `From {IF}
  to {BI}', \emph{Synthese}, 167 (2009), 2, 207--230.

\bibitem{vanDalen01}
\textsc{van Dalen, D.}, `Intuitionistic logic', in L.~Goble, (ed.), \emph{The
  Blackwell Guide to Philosophical Logic}, Blackwell, 2001, pp. 224--257.

\bibitem{Enderton1970}
\textsc{Enderton, H.B.}, `Finite partially-ordered quantifiers',
  \emph{Zeitschrift fur Mathematische Logik und Grundlagen der Mathematik},
  (1970), 16, 393--397.

\bibitem{HintikkaSandu1989}
\textsc{Hintikka, J.}, and \textsc{G.~Sandu}, `Informational independence as a
  semantical phenomenon', in J.~E. Fenstad, I.~T. Frolov, and R.~Hilpinen,
  (eds.), \emph{Logic, Methodology and Philosophy of Science}, vol.~8,
  Amsterdam: Elsevier, 1989, pp. 571--589.

\bibitem{Hodges1997a}
\textsc{Hodges, W.}, `Compositional semantics for a langauge of imperfect
  information', \emph{Logic Journal of the IGPL}, 5 (1997), 539--563.

\bibitem{Hodges1997b}
\textsc{Hodges, W.}, `Some strange quantifiers', in J.~Mycielski, G.~Rozenberg,
  and A.~Salomaa, (eds.), \emph{Structures in Logic and Computer Science: A
  Selection of Essays in Honor of A. Ehrenfeucht}, vol. 1261 of \emph{Lecture
  Notes in Computer Science}, London: Springer, 1997, pp. 51--65.

\bibitem{KontinenNurmi2009}
\textsc{Kontinen, J.}, and \textsc{V. Nurmi}, `Team logic and second-order
  logic', \emph{Fundamenta Informaticae}, 106 (2011), 259--272.

\bibitem{KontVan09}
\textsc{Kontinen, J.}, and \textsc{J. V\"{a}\"{a}n\"{a}nen}, `On
  definability in dependence logic', \emph{Journal of Logic, Language and
  Information}, 18(3) (2009), 317--332.

\bibitem{Ville_thesis}
\textsc{Nurmi, V.}, \emph{Dependence Logic: Investigations into Higher-Order
  Semantics Defined on Teams}, Ph.D. thesis, University of Helsinki, 2009.

\bibitem{OHearnPym1999}
\textsc{O'Hearn, P.}, and \textsc{D.~Pym}, `The logic of bunched implications',
  \emph{Bulletin of Symbolic Logic}, 5(2) (1999), 215--244.

\bibitem{Pym2002}
\textsc{Pym, D.}, \emph{The Semantics and Proof Theory of the Logic of Bunched
  Implications}, Kluwer Academic Publishers, 2002.

\bibitem{Van07dl}
\textsc{V{\"a}{\"a}n{\"a}nen, J.}, \emph{Dependence Logic: A New Approach to
  Independence Friendly Logic}, Cambridge: Cambridge University Press, 2007.

\bibitem{Walkoe1970}
\textsc{Walkoe, W.J.}, `Finite partially-ordered quantification', \emph{Journal
  of Symbolic Logic}, 35 (1970), 535--555.

\bibitem{Yang2010}
\textsc{Yang, F.}, \emph{On Definability in Intuitionistic Dependence Logic},
  manuscript, 2010.


\end{thebibliography}
\end{document}